\documentclass{amsart}
\theoremstyle{plain}
\newtheorem{Thm}{Theorem}

\newtheorem{Cor}[Thm]{Corollary}

\renewcommand{\theMain}{}

\errorcontextlines=0 \numberwithin{equation}{section}

\newcommand{\rad}{\operatorname{rad}}
\newcommand{\ep}{\epsilon}
\newcommand{\bt}{\bigtriangleup}
\newcommand{\btd}{\bigtriangledown}
\newcommand{\Om}{\Omega}

\begin{document}

%begin Topmatter
\title[non-local heat flows]
{non-local heat flows and gradient estimates on closed manifolds}

\author{Li Ma, Liang Cheng}

\address{Department of mathematical sciences \\
Tsinghua university \\
Beijing 100084 \\
China} \email{lma@math.tsinghua.edu.cn} \dedicatory{}
\date{May 26th, 2009}

\begin{abstract}

In this paper, we study two kind of $L^2$ norm preserved non-local
heat flows on closed manifolds. We first study the global existence,
stability and asymptotic behavior to such non-local heat flows. Next
we give the gradient estimates of positive solutions to these heat
flows.


\end{abstract}

\thanks{$^*$ The research is partially supported by the National Natural Science
Foundation of China 10631020 and SRFDP 20060003002. }
 \maketitle

\section{Introduction}
In this paper, we consider the global existence, stability,
asymptotic behavior and gradient estimates for two kind of $L^2$
preserving heat flow which have positive solutions on closed
manifolds. This is a continuation of our earlier study of
non-local heat flows in \cite{MC} and \cite{MA}. Our work is also
motivated by the recent work of C.Caffarelli and F.Lin
\cite{CL09}, where they have studied the global existence and
regularity of $L^2$ norm preserving heat flow such as $$\partial_t
u=\Delta u+\lambda(t)u$$ with
$$\lambda(t)=\frac{\int_{\Omega}|\nabla
u|^2dx}{\int_{\Omega}u^2dx}.$$ They also extend the method to
study a family of singularly perturbed systems of non-local
parabolic equations. We remark that the non-local heat flow
naturally arises in geometry such that the flow preserves some
$L^p$ norm in the sense that some the geometrical quantity (such
as length, area and so on) is preserved in the geometric heat
flow. For more references on geometric flows such as harmonic map
heat flows and non-local heat flows, one may see \cite{A98},
\cite{Str},\cite{MA} and \cite{MC}.

We firstly study the following linear heat flow on a closed
Riemannian manifold $M$,
\begin{equation*}
\left\{
\begin{array}{l}
         \partial_t u=\Delta u+\lambda(t)u+A(x,t) \ \ \ \text{in}\ M\times\mathbb{R}_{+}, \\
          u(x,0)=g(x) \ \ \ \text{in}\ M,
\end{array}
\right.
\end{equation*}
where $g\in C^1(M)$, $A(x,t)$ is a given non-negative smooth
function, and $\lambda(t)$ is chosen such that the flow preserves
the $L^2$ the norm of the solution. In fact, a direct computation
shows that
$$\frac{1}{2}\frac{d}{dt}\int_{M}u^2dx=\int_{M}u u_t=-\int_{M}|\nabla u|^2dx+\lambda(t)\int_{M}u^2dx+\int_{M}uAdx.$$
Thus, one has $\lambda(t)=\frac{\int_{M}(|\nabla
u|^2-uA)dx}{\int_{M}g^2dx}$ to preserve the $L^2$ norm. Without loss
of generality we may assume $\int_{M}g^2dx=1$. Thus we are lead to
consider the following problem on the closed Riemannian manifold $M$
\begin{equation} \label{eq1}
\left\{
\begin{array}{l}
         \partial_t u=\Delta u+\lambda(t)u+A(x,t) \ \ \ \text{in}\ M\times\mathbb{R}_{+}, \\
          u(x,0)=g(x) \ \ \ \text{in}\ M,
\end{array}
\right.
\end{equation}
where $\lambda(t)=\int_{M}(|\nabla u|^2-uA)dx$,  $g(x)\geq 0$,
$\int_{M}g^2dx=1$.

We next study the following nonlinear heat flow on the closed
Riemannian manifold $M$,
\begin{equation*}
\left\{
\begin{array}{l}
         \partial_t u=\Delta u+\lambda(t)u-u^p \ \ \ \text{in}\ M\times\mathbb{R}_{+}, \\
          u(x,0)=g(x) \ \ \ \text{in}\ M,
\end{array}
\right.
\end{equation*}
where $p>1$, which has the positive solution and preserves the
$L^2$ the norm. Likewise,
$$\frac{1}{2}\frac{d}{dt}\int_{M}u^2dx=\int_{M}u
u_t=-\int_{M}|\nabla
u|^2dx+\lambda(t)\int_{M}u^2dx-\int_{M}u^{p+1}dx.$$ Thus, one must
have $\lambda(t)=\frac{\int_{M}(|\nabla
u|^2+u^{p+1})dx}{\int_{M}g^2dx}$ to preserve the $L^2$ norm. Without
loss of generality we assume $\int_{M}g^2dx=1$. Then we consider the
following problem on closed Riemannian manifold $M$
\begin{equation} \label{eq2}
\left\{
\begin{array}{l}
         \partial_t u=\Delta u+\lambda(t)u-u^p\ \ \ \text{in}\ M\times\mathbb{R}_{+} \\
          u(x,0)=g(x) \ \ \ \text{in}\ M
\end{array}
\right.
\end{equation}
where $p>1$, $\lambda(t)=\int_{M}(|\nabla u|^2+u^{p+1}) dx$,
$g(x)\geq 0\ \text{in}\ M$, $\int_{M}g^2dx=1$ and $g\in C^1(M)$.

Similar to the global existence results obtained in C.Caffarelli
and F.Lin \cite{CL09}, we have following two results.
\begin{Thm}\label{thm1}
Problem (\ref{eq1}) has a global non-negative solution $u(t)\in
L^{\infty}(\mathbb{R}_{+},H^1(M))\cap
L_{loc}^2(\mathbb{R}_{+},H^2(M))$ if $A\in
L^{\infty}(\mathbb{R}_{+},H^1(M))$ .
\end{Thm}

\begin{Thm}\label{thm2}
Problem (\ref{eq2}) has a global positive solution
$$u(t)\in
L^{\infty}(\mathbb{R}_{+},H^1(M))\cap
L^{\infty}(\mathbb{R}_{+},L^{p+1}(M))\cap
L_{loc}^2(\mathbb{R}_{+},H^2(M)).$$
\end{Thm}

We also have the stability results for both problem (\ref{eq1})
and (\ref{eq2}).

\begin{Thm}\label{thm5}
Let $u,v$ be the two non-negative solutions to problem (\ref{eq1})
with initial data $g_u,g_v$ at $t=0$. Suppose $A\in
L^{\infty}(\mathbb{R}_{+},H^1(M))$. Then
$$||u-v||_{L^2}^2\leq
||g_u-g_v||_{L^2}^2 \exp(C_1 t)$$ and $$||\nabla(u-v)||_{L^2}^2\leq
||\nabla(g_u-g_v)||_{L^2}^2 \exp(C_2 t),$$
 where $C_i$, $i=1,2$, are constants
depending on the upper bound of $||g_u||_{H^1(M)}, ||g_v||_{H^1(M)}$
and $||A||_{L^{\infty}(\mathbb{R}_{+},H^1(M))}$. In particular, the
solution to problem (\ref{eq1}) is unique.
\end{Thm}

\begin{Thm}\label{thm6}
Let $u,v$ be the two bounded positive solutions to problem
(\ref{eq2}) with initial data $g_u,g_v$ at $t=0$, where $g_u,g_v\in
H^1(M)\cap L^{\infty}(M)$. Then $$||u-v||_{L^2}^2\leq
||g_u-g_v||_{L^2}^2 \exp(C_1 t)$$ and $$||\nabla(u-v)||_{L^2}^2\leq
||\nabla(g_u-g_v)||_{L^2}^2 \exp(C_2 t),$$ where $C_i$, $i=1,2$, are
constants depending on the upper bound of $||g_u||_{H^1(M)},
||g_v||_{H^1(M)}$ and $||g_u||_{L^{\infty}}, ||g_v||_{L^{\infty}}$.
In particular, the solution to problem (\ref{eq2}) is unique.
\end{Thm}

As the simple applications to Theorem \ref{thm1} and theorem
\ref{thm2}, we can study asymptotic behavior of $u(t)$ of problem
(\ref{eq1}) and problem (\ref{eq2}). We have the following two
corollaries.
\begin{Cor}\label{cor1}
Suppose $A\in L^{\infty}(\mathbb{R}_{+},H^1(M))\cap
L^{2}(\mathbb{R}_{+},L^2(M))$ in Theorem \ref{thm1}. Let $u(t)$ be
the solution to problem (\ref{eq1}). Then one can take $t_i\to
\infty$ such that $\lambda(t_i)\to\lambda_{\infty}$,
$u(x,t_i)\rightharpoonup u_{\infty}(x)$ in $H^1(M)$ and $u_{\infty}$
solves the equation $\Delta
u_{\infty}+\lambda_{\infty}u_{\infty}+A=0$ in $M$ with
$\int_{M}|u_{\infty}|^2dx=1$.
\end{Cor}

\begin{Cor}\label{cor2}
Suppose $u(t)$ is the positive solution to problem (\ref{eq2}). Then
one can take $t_i\to \infty$ such that
$\lambda(t_i)\to\lambda_{\infty}$, $u(x,t_i)\rightharpoonup
u_{\infty}(x)$ in $H^1(M)$ and $u_{\infty}>0$ solves the equation
$\Delta u_{\infty}+\lambda_{\infty}u_{\infty}-u_{\infty}^p=0$ in $M$
with $\int_{M}|u_{\infty}|^2dx=1$.
\end{Cor}

In the study of geometric analysis as well as other elliptic or
parabolic equations, it is well known that the Harnack inequality
plays a important role (see for instance \cite{BH97}, \cite{P02}
and \cite{SY}). As showed in \cite{SY}, the Harnack inequality for
positive solutions is a consequence of the gradient estimates for
positive solutions. We have the following gradient estimates for
the type of heat equations related to problem (\ref{eq1}) and
(\ref{eq2}).
\begin{Thm}\label{thm3}
Suppose $M$ is a closed Riemannian manifold with Ricci curvature
bounded from below by $\geq -K$. Let $u>0$ be a smooth solution to
the heat equation on $M\times [0,T)$
$$
(\partial_t-\Delta) u=\lambda(t)u+A(x,t).
$$
Let, for $a>1$ and $w=\log u$,
$$
F=t(|\nabla w|^2-aw_t+a(\lambda+u^{-1}A)).
$$
Then there is a constant $C(u^{-1},|A|,|\nabla A|,|\Delta A|,K,a,T)$
such that
$$\sup_{M\times[0,T]}F\leq C(|u|^{-1},|A|,|\nabla A|,|\Delta A|,K,a,T).$$
\end{Thm}

\begin{Thm}\label{thm4}
Suppose $M$ is a closed Riemannian manifold with Ricci curvature
bounded from below by $\geq -K$. Let $u>0$ be a smooth solution to
the heat equation on $M\times [0,T)$
$$
(\partial_t-\Delta) u=\lambda(t)u-u^p.
$$
Let, for $a>1$ and $w=\log u$, $$ F=t(|\nabla
w|^2-aw_t+a(\lambda-u^{p-1})).
$$ Then there is a constant  $C(u^{p-1},K,a,p,T)$
such that
$$\sup_{M\times[0,T]}F\leq C(|u|^{p-1},K,a,p,T).$$
\end{Thm}

The proofs of the results above are similar to our earlier work
\cite{LM}.

This paper is organized as follows. In section \ref{sect2} we
study the global existence, stability and asymptotic behavior of
solutions to the problem (\ref{eq1}) and problem (\ref{eq2}). In
particular, we give the proofs of theorem \ref{thm1} to theorem
\ref{thm6}, corollary \ref{cor1} and corollary \ref{cor2}. In
section \ref{sect3} we do the gradient estimates for problem
(\ref{eq1}) and problem (\ref{eq2}).

\section{global existence and stability  property}\label{sect2}
In this section we study the global existence, stability and
asymptotic behavior of solutions to the problem (\ref{eq1}) and
problem (\ref{eq2}). For global existence, we use the idea of
\cite{CL09} theorem 1.1 and construct a series of solutions to
linear parabolic equations to converge to the solution of the
non-local heat flow.

\textbf{Proof of Theorem \ref{thm1}}. Let us define a series
non-negative functions $u^{(k)}$ by
\begin{eqnarray} \label{eq3}
\left\{
\begin{array}{l}
         u^{(0)}=g, \lambda^{(k)}(t)=\int_{M}(|\nabla u^{(k)}|^2-u^{(k)}A) dx, \\
          \partial_t u^{(k+1)}=\Delta u^{(k+1)}+\lambda^{(k)}(t)u^{(k+1)}+A ,        \\
          u^{(k+1)}(x,0)=g(x),
\end{array}
\right.
\end{eqnarray}
which are a series of initial boundary value problems of linear
parabolic systems.

To prove the convergence of the series $\{u^{(k)}\}$ constructed
above, we estimate for $k\geq 0$
\begin{equation}\label{aaa}
\frac{1}{2}\frac{d}{dt}\int_{M}|u^{(k+1)}|^2 dx+\int_{M}|\nabla
u^{(k+1)}|^2 dx= \lambda^{(k)}(t)\int_{M}|u^{(k+1)}|^2 dx
\end{equation}
$$
+\int_{M}u^{(k+1)} Adx,
$$
\begin{equation}\label{bbb}
\frac{1}{2}\frac{d}{dt}\int_{M}|\nabla u^{(k+1)}|^2
dx+\int_{M}|\Delta u^{(k+1)}|^2 dx=\lambda^{(k)}(t)\int_{M}|\nabla
u^{(k+1)}|^2 dx
\end{equation}
$$
+\int_{M}\nabla u^{(k+1)}\cdot \nabla Adx,
$$
\begin{equation}\label{ccc}
\frac{1}{2}\frac{d}{dt}\int_{M}|\nabla u^{(k+1)}|^2
dx+\int_{M}|u_t^{(k+1)}|^2
dx=\frac{\lambda^{(k)}(t)}{2}\frac{d}{dt}\int_{M}|u^{(k+1)}|^2 dx
\end{equation}
$$
+\int_{M}u_t^{(k+1)}Adx.
$$
Since $\lambda^{(k)}(t)\leq \int_{M}|\nabla u^{(k)}|^2$, by
(\ref{bbb}), we have
\begin{eqnarray*}
\frac{1}{2}\frac{d}{dt}\int_{M}|\nabla u^{(k+1)}|^2
dx&\leq&\lambda^{(k)}(t)\int_{M}|\nabla u^{(k+1)}|^2 dx
+\int_{M}\nabla u^{(k+1)}\cdot \nabla Adx\\
&\leq&(\int_{M}|\nabla u^{(k)}|^2)\int_{M}|\nabla u^{(k+1)}|^2 dx
+\frac{1}{2}\int_{M}|\nabla u^{(k+1)}|^2dx\\
&+&\frac{1}{2}\int_{M}|\nabla A|^2dx.
\end{eqnarray*}
Now we define
$c_1=max(1,||A||^2_{L^{\infty}(\mathbb{R}_{+},H^1(M))})$, we get
$$
\frac{d}{dt}(\int_{M}|\nabla u^{(k+1)}|^2dx+c_1)\leq
2(\int_{M}|\nabla u^{(k)}|^2dx+c_1)(\int_{M}|\nabla
u^{(k+1)}|^2dx+c_1).
$$
Hence
$$
\int_{M}|\nabla u^{(k+1)}|^2dx+c_1\leq (\int_{M}|\nabla
g|^2dx+c_1)\exp(2\int^{t}_0(\int_{M}|\nabla u^{(k)}|^2dx+c_1)dt).
$$
By induction, there is $\delta$ depending only on $\int_{M}|\nabla
g|^2dx$ and $||A||_{L^{\infty}(\mathbb{R}_{+},H^1(M))}$ such that
\begin{equation}\label{estimate1}
\int_{M}|\nabla u^{(k+1)}|^2dx\leq c'_1,\ \ \text{for}\
t\in[0,\delta], k\geq 1,
\end{equation}
where $c'_1$ is a constant depending on $\int_{M}|\nabla g|^2dx$ and
$||A||_{L^{\infty}(\mathbb{R}_{+},H^1(M))}$. Hence
$\lambda^{(k+1)}\leq \int_{M}|\nabla u^{(k+1)}|^2dx\leq c'_1$. By
(\ref{aaa}), we have
\begin{eqnarray*}
\frac{d}{dt}\int_{M}|u^{(k+1)}|^2 dx&\leq&
2\lambda^{(k)}(t)\int_{M}|u^{(k+1)}|^2 dx+2\int_{M}u^{(k+1)} Adx\\
&\leq& (2\lambda^{(k)}(t)+1)\int_{M}|u^{(k+1)}|^2
dx+\int_{M}|A|^2dx.
\end{eqnarray*}
Hence,
\begin{equation}\label{estimate2}
\int_{M}|u^{(k+1)}|^2dx\leq c_2,\ \ \text{for}\ t\in[0,\delta],
k\geq1,
\end{equation}
where $c_2$ depending on $\int_{M}|g|^2dx$ and
$||A||_{L^{\infty}(\mathbb{R}_{+},H^1(M))}$. Now integrate
(\ref{bbb}) with t, we get
$$
\frac{1}{2}\int_{M}|\nabla u(t)^{(k+1)}|^2
dx-\frac{1}{2}\int_{M}|\nabla u(0)^{(k+1)}|^2
dx+\int^{\delta}_0\int_{M}|\Delta u^{(k+1)}|^2 dxdt
$$
$$
=\int^{\delta}_0\lambda^{(k)}(t)\int_{M}|\nabla u^{(k+1)}|^2 dxdt
+\int^{\delta}_0\int_{M}\nabla u^{(k+1)}\cdot \nabla Adxdt.
$$
We have
$$
\int^{\delta}_0\int_{M}|\Delta u^{(k+1)}|^2 dxdt
\leq\frac{1}{2}\int_{M}|\nabla g|^2
dx+\int^{\delta}_0\lambda^{(k)}(t)\int_{M}|\nabla u^{(k+1)}|^2 dxdt
$$
$$
+\int^{\delta}_0\int_{M}\nabla u^{(k+1)}\cdot \nabla Adxdt.
$$
Hence,
\begin{equation}\label{estimate3}
\int^{\delta}_0\int_{M}|\Delta u^{(k+1)}|^2 dxdt\leq c_3,
\end{equation}
where $c_3$ depending on $\int_{M}|\nabla g|^2dx$ and
$||A||_{L^{\infty}(\mathbb{R}_{+},H^1(M))}$. We integrate
(\ref{ccc}) with t and we get
$$
\frac{1}{2}\int_{M}|\nabla u(t)^{(k+1)}|^2
dx-\frac{1}{2}\int_{M}|\nabla u(0)^{(k+1)}|^2
dx+\int^{\delta}_0\int_{M}|u_t^{(k+1)}|^2 dxdt
$$
$$
=\int^{\delta}_0\frac{\lambda^{(k)}(t)}{2}\frac{d}{dt}\int_{M}|u^{(k+1)}|^2
dxdt+\int^{\delta}_0\int_{M}u_t^{(k+1)}Adxdt.
$$
We then have
$$
\int^{\delta}_0\int_{M}|u_t^{(k+1)}|^2 dxdt \leq
\frac{1}{2}\int_{M}|\nabla g|^2 dx+
\int^{\delta}_0\frac{\lambda^{(k)}(t)}{2}\frac{d}{dt}\int_{M}|u^{(k+1)}|^2
dxdt
$$
$$
+\int^{\delta}_0\int_{M}u_t^{(k+1)}Adxdt.
$$
Hence
\begin{equation}\label{estimate4}
\int^{\delta}_0\int_{M}|u_t^{(k+1)}|^2 dxdt\leq c_4,
\end{equation}
where $c_4$ depending on $\int_{M}|\nabla g|^2dx$ and
$||A||_{L^{\infty}(\mathbb{R}_{+},H^1(M))}$.

By (\ref{estimate1}), (\ref{estimate2}), (\ref{estimate3}) and
(\ref{estimate4}), there is a subsequence of $\{u^{(k)}\}$ (still
denoted by $\{u^{(k)}\}$) and a function $u(t)\in
L^{\infty}([0,\delta],H^1(M))\cap L^2([0,\delta],H^2(M))$ with
$\partial_t u(t)\in L^2([0,\delta],L^2(M))$ such that
$u^{(k)}\rightharpoonup u$ weak$^\ast$ in
$L^{\infty}([0,\delta],H^1(M))$ and weakly in
$L^2([0,\delta],H^2(M))$. Then we have $u^{(k)}\to u$ strongly in
$L^2([0,\delta],H^1(M))$ and $u(t)\in C([0,\delta],L^2(M))$. Hence
$\lambda^{(k)}(t)\to \lambda(t)$ strongly in $L^2([0,\delta])$.
Thus, we get a local strong solution to problem (\ref{eq1}). Next,
starting from $t=\delta$ we can extend the local solution to
$[0,2\delta]$ in exactly the same way as above. By induction, we
have a global solution to problem (\ref{eq1}).$\Box$

The proof of Theorem \ref{thm2} is similar to Theorem \ref{thm1}
with slight modification.

\textbf{Proof of Theorem \ref{thm2}.} Let us define a series
$u^{(k)}$ by
\begin{eqnarray} \label{eq4}
\left\{
\begin{array}{l}
         u^{(0)}=g, \lambda^{(k)}(t)=\int_{M}(|\nabla u^{(k)}|^2+(u^{(k)})^{p+1}) dx, \\
          \partial_t u^{(k+1)}=\Delta u^{(k+1)}+ \lambda^{(k)}(t)u^{(k+1)}-(u^{k+1})^p,        \\
          u^{(k+1)}(x,0)=g(x),
\end{array}
\right.
\end{eqnarray}
which are a series of initial boundary value problems of linear
parabolic systems. By the maximum principle, we know that
$u^{(k)}>0$.

To prove the convergence of series $\{u^{(k)}\}$ constructed above,
we estimate for $k\geq 0$
\begin{equation}\label{aaa1}
\frac{1}{2}\frac{d}{dt}\int_{M}|u^{(k+1)}|^2 dx+\int_{M}(|\nabla
u^{(k+1)}|^2+(u^{(k+1)})^{p+1})dx
\end{equation}
$$
= \lambda^{(k)}(t)\int_{M}|u^{(k+1)}|^2 dx,
$$
\begin{equation}\label{bbb1}
\frac{1}{2}\frac{d}{dt}\int_{M}|\nabla u^{(k+1)}|^2
dx+\int_{M}|\Delta u^{(k+1)}|^2 dx+\int_{M}p(u^{(k+1)})^{p-1}|\nabla
u^{(k+1)}|^2dx
\end{equation}
$$
=\lambda^{(k)}(t)\int_{M}|\nabla u^{(k+1)}|^2 dx,
$$
\begin{equation}\label{ccc1}
\frac{1}{2}\frac{d}{dt}\int_{M}|\nabla u^{(k+1)}|^2
dx+\int_{M}|u_t^{(k+1)}|^2
dx+\frac{1}{p+1}\frac{d}{dt}\int_{M}(u^{(k+1)})^{p+1}dx
\end{equation}
$$
=\frac{\lambda^{(k)}(t)}{2}\frac{d}{dt}\int_{M}|u^{(k+1)}|^2 dx.
$$
\begin{equation}\label{ddd1}
\frac{1}{p+1}\frac{d}{dt}\int_{M}(u^{(k+1)})^{p+1}
dx+\int_{M}p(u^{k+1})^{p-1}|\nabla
u^{(k+1)}|^2dx+\int_{M}(u^{(k+1)})^{2p}dx
\end{equation}
$$
=\lambda^{(k)}(t)\int_{M}(u^{(k+1)})^{p+1}dx.
$$
By (\ref{bbb1}) and (\ref{ddd1}), we have
$\frac{d}{dt}\int_{M}|\nabla u^{(k+1)}|^2 dx\leq
2\lambda^{(k)}(t)\int_{M}|\nabla u^{(k+1)}|^2 dx$ and
$\frac{d}{dt}\int_{M}(u^{(k+1)})^{p+1} dx \leq (p+1)
\lambda^{(k)}(t)\int_{M}(u^{(k+1)})^{p+1}dx$. Hence
$\frac{d}{dt}\lambda^{(k+1)}\leq (p+1)\lambda^{(k)}\lambda^{(k+1)}$
and $\lambda^{(k+1)}\leq \int_{M}(|\nabla g|^2+g^{p+1}) dx\cdot
\exp((p+1)\int^{t}_0\lambda^{(k+1)}dt)$. By induction, there is
$\delta$ depending only on $\int_{M}|\nabla g|^2dx$ and
$\int_{M}g^{p+1}dx$ such that
\begin{equation}\label{estimate11}
\lambda^{(k+1)}(t)\leq c_5,\ \ \text{for}\ t\in[0,\delta], k\geq1,
\end{equation}
where $c_5$ is a constant depending on $\int_{M}|\nabla g|^2dx$ and
$\int_{M}g^{p+1}dx$. Hence
\begin{equation}\label{estimate12}
\int_{M}|\nabla u^{(k+1)}|^2dx\leq c_5,\ \int_{M}|
u^{(k+1)}|^{p+1}dx\leq c_5 \ \ \text{for}\ t\in[0,\delta], k\geq1.
\end{equation}
Now we integrate (\ref{bbb1}) with t and we get
$$
\frac{1}{2}\int_{M}|\nabla u(t)^{(k+1)}|^2
dx-\frac{1}{2}\int_{M}|\nabla u(0)^{(k+1)}|^2
dx+\int^{\delta}_0\int_{M}|\Delta u^{(k+1)}|^2 dxdt
$$
$$
+\int^{\delta}_0\int_{M}p(u^{(k+1)})^{p-1}|\nabla
u^{(k+1)}|^2dxdt=\int^{\delta}_0\lambda^{(k)}(t)\int_{M}|\nabla
u^{(k+1)}|^2 dxdt.
$$
We then have
$$
\int^{\delta}_0\int_{M}|\Delta u^{(k+1)}|^2 dxdt
\leq\frac{1}{2}\int_{M}|\nabla g|^2
dx+\int^{\delta}_0\lambda^{(k)}(t)\int_{M}|\nabla u^{(k+1)}|^2 dxdt.
$$
Hence
\begin{equation}\label{estimate13}
\int^{\delta}_0\int_{M}|\Delta u^{(k+1)}|^2 dxdt\leq c_6,
\end{equation}
where $c_6$ depending on $\int_{M}|g|^2dx$ and $\int_{M}|\nabla g|^2
dx$. Integrating (\ref{ccc1}) with t, we get
$$
\frac{1}{2}\int_{M}|\nabla u(t)^{(k+1)}|^2
dx-\frac{1}{2}\int_{M}|\nabla u(0)^{(k+1)}|^2
dx+\int^{\delta}_0\int_{M}|u_t^{(k+1)}|^2 dxdt
$$
$$
+\frac{1}{p+1}\int_{M}(u^{(k+1)}(t))^{p+1}dx-\frac{1}{p+1}\int_{M}g^{p+1}dx=
\int^{\delta}_0\frac{\lambda^{(k)}(t)}{2}\frac{d}{dt}\int_{M}|u^{(k+1)}|^2
dxdt.
$$
We then have
$$
\int^{\delta}_0\int_{M}|u_t^{(k+1)}|^2 dxdt \leq
\frac{1}{2}\int_{M}|\nabla g|^2 dx+\frac{1}{p+1}\int_{M}g^{p+1}dx+
\int^{\delta}_0\frac{\lambda^{(k)}(t)}{2}\frac{d}{dt}\int_{M}|u^{(k+1)}|^2
dxdt.
$$
Hence
\begin{equation}\label{estimate14}
\int^{\delta}_0\int_{M}|u_t^{(k+1)}|^2 dxdt\leq c_7,
\end{equation}
where $c_7$ depending on $\int_{M}|g|^2dx$, $\int_{M}|\nabla g|^2dx$
and $\int_{M}g^{p+1}dx$.

By (\ref{estimate11}), (\ref{estimate12}), (\ref{estimate13}) and
(\ref{estimate14}), there is a subsequence of $\{u^{(k)}\}$ (still
denoted by $\{u^{(k)}\}$) and a function $u(t)\in
L^{\infty}([0,\delta],H^1(M))\cap L^2([0,\delta],H^2(M))\cap
L^{\infty}([0,\delta],L^{p+1}(M))$ with $\partial_t u(t)\in
L^2([0,\delta],L^2(M))$ such that $u^{(k)}\rightharpoonup u$
weak$^\ast$ in $L^{\infty}([0,\delta],H^1(M))$, weakly in
$L^2([0,\delta],H^2(M))$ and weakly in
$L^{\infty}([0,\delta],L^{p+1}(M))$. Then we have $u^{(k)}\to u$
strongly in $L^2([0,\delta],H^1(M))$ and $u(t)\in
C([0,\delta],L^2(M))$. Hence $\lambda^{(k)}(t)\to \lambda(t)$
strongly in $L^2([0,\delta])$. Thus, we get a local strong solution
to problem (\ref{eq2}). Next, starting from $t=\delta$ we can extend
the local solution to $[0,2\delta]$ in exactly the same way as
above. By induction, we have a global solution to problem
(\ref{eq2}).
 $\Box$

Now we can study the stability for the problem ({\ref{eq1}) and
problem (\ref{eq2}).

\textbf{Proof of Theorem \ref{thm5}.} By the arguments in Theorem
\ref{thm1}, we can take a constant $C$ such that all
$||u||_{L^{\infty}(\mathbb{R}_{+},H^1(M))}$,
$||v||_{L^{\infty}(\mathbb{R}_{+},H^1(M))}$,
$||A||_{L^{\infty}(\mathbb{R}_{+},H^1(M))}$,
$||\lambda_u(t)||_{L^{\infty}(\mathbb{R}_{+})}$ and
$||\lambda_v(t)||_{L^{\infty}(\mathbb{R}_{+})}$ not less than $C$,
where $C$ is only depending on upper bound of $||g_u||_{H^1(M)},
||g_v||_{H^1(M)}$ and $||A||_{L^{\infty}(\mathbb{R}_{+},H^1(M))}$.
First we calculate
\begin{eqnarray*}
\frac{1}{2}\frac{d}{dt}\int_{M}(u-v)^2dx
&=&\int_{M}(u-v)(u_t-v_t)dx\\
&=&\int_{M}(u-v)(\Delta(u-v)+\lambda_u(t)u-\lambda_v(t)v)dx\\
&=&-\int_{M}|\nabla(u-v)|^2dx+\int_{M}(u-v)(\lambda_u(t)u-\lambda_v(t)v)dx
\end{eqnarray*}
Note that
\begin{eqnarray*}
&&\int_{M}(u-v)(\lambda_u(t)u-\lambda_v(t)v)dx\\
&=&(\lambda_u(t)-\lambda_v(t))\int_{M}(u-v)udx+\lambda_v(t)\int_{M}(u-v)^2dx\\
&\leq&|\lambda_u(t)-\lambda_v(t)|(\int_{M}(u-v)^2dx)^{\frac{1}{2}}(\int_{M}u^2dx)^{\frac{1}{2}}
+\lambda_v(t)\int_{M}(u-v)^2dx\\
&\leq&C|\lambda_u(t)-\lambda_v(t)|(\int_{M}(u-v)^2dx)^{\frac{1}{2}}
+C\int_{M}(u-v)^2dx,
\end{eqnarray*}
and
\begin{eqnarray}\label{lambda}
|\lambda_u(t)-\lambda_v(t)|
&=&|\int_{M}((|\nabla u|^2-|\nabla v|^2)-(uA-vA))dx|\nonumber\\
&\leq&\int_{M}|\nabla(u-v)|(|\nabla u|+|\nabla
v|)dx+|\int_{M}(u-v)Adx|\nonumber\\
&\leq&(\int_{M}|\nabla(u-v)|^2dx)^{\frac{1}{2}}(\int_{M}(|\nabla
u|+|\nabla
v|)^2dx)^{\frac{1}{2}}\nonumber\\
&&+(\int_{M}(u-v)^2dx)^{\frac{1}{2}}(\int_{M}A^2dx)^{\frac{1}{2}}\nonumber\\
&\leq&C(\int_{M}|\nabla(u-v)|^2dx)^{\frac{1}{2}}+C(\int_{M}(u-v)^2dx)^{\frac{1}{2}}.
\end{eqnarray}
We have
\begin{eqnarray*}
&&\frac{1}{2}\frac{d}{dt}\int_{M}(u-v)^2dx\\
&\leq&-\int_{M}|\nabla(u-v)|^2dx+C^2(\int_{M}|\nabla(u-v)|^2dx)^{\frac{1}{2}}(\int_{M}(u-v)^2dx)^{\frac{1}{2}}\\
&&+C^2\int_{M}(u-v)^2dx+C\int_{M}(u-v)^2dx\\
&\leq&-\frac{1}{2}\int_{M}|\nabla(u-v)|^2dx+(\frac{C^4}{2}+C^2+C)\int_{M}(u-v)^2dx.
\end{eqnarray*}
By Gronwall inequality, we have
$$
||u-v||_{L^2}^2\leq ||g_u-g_v||_{L^2}^2
\exp((\frac{C^4}{2}+C^2+C)t).
$$
Further more,
\begin{eqnarray*}
\frac{1}{2}\frac{d}{dt}\int_{M}|\nabla(u-v)|^2dx
&=&\int_{M}\nabla(u-v) \cdot \nabla(u-v)_tdx\\
&=&-\int_{M}\Delta(u-v) \cdot (u-v)_tdx\\
&=&-\int_{M}\Delta(u-v)\cdot(\Delta(u-v)+\lambda_u(t)u-\lambda_v(t)v)dx\\
&=&-\int_{M}(\Delta(u-v))^2dx+\int_{M}\nabla(u-v)\cdot\nabla(\lambda_u(t)u-\lambda_v(t)v)dx
\end{eqnarray*}
Likewise,
\begin{eqnarray*}
&&\int_{M}\nabla(u-v)\cdot\nabla(\lambda_u(t)u-\lambda_v(t)v)dx\\
&=&(\lambda_u(t)-\lambda_v(t))\int_{M}\nabla(u-v)\cdot \nabla udx+\lambda_v(t)\int_{M}|\nabla(u-v)|^2dx\\
&\leq&|\lambda_u(t)-\lambda_v(t)|(\int_{M}|\nabla(u-v)|^2dx)^{\frac{1}{2}}(\int_{M}|\nabla
u|^2dx)^{\frac{1}{2}}
+\lambda_v(t)\int_{M}|\nabla(u-v)|^2dx\\
&\leq&C|\lambda_u(t)-\lambda_v(t)|(\int_{M}|\nabla(u-v)|^2dx)^{\frac{1}{2}}
+C\int_{M}|\nabla(u-v)|^2dx.
\end{eqnarray*}
By Poincare inequality and (\ref{lambda}), we have
\begin{eqnarray*}
|\lambda_u(t)-\lambda_v(t)|
&\leq&C'(\int_{M}|\nabla(u-v)|^2dx)^{\frac{1}{2}},
\end{eqnarray*}
where $C'$ is a constant depending on the constant of Poincare
inequality and $C$. We take a constant of the maximum of $C$ and
$C'$, and still denote it $C$ for convenience. Then we have
\begin{eqnarray*}
\frac{1}{2}\frac{d}{dt}\int_{M}|\nabla(u-v)|^2dx \leq
(C^2+C)\int_{M}|\nabla(u-v)|^2dx.
\end{eqnarray*}
By Gronwall inequality, we have
$$
||\nabla(u-v)||_{L^2}^2\leq ||\nabla(g_u-g_v)||_{L^2}^2
\exp((C^2+C)t).
$$
 $\Box$

\textbf{Proof of Theorem \ref{thm6}.} By the arguments in Theorem
\ref{thm2}, we can take a constant $C$ such that all
$||u||_{L^{\infty}(\mathbb{R}_{+},H^1(M))}$,
$||v||_{L^{\infty}(\mathbb{R}_{+},H^1(M))}$,
$||u||_{L^{\infty}(\mathbb{R}_{+},L^{\infty}(M))}$,
$||v||_{L^{\infty}(\mathbb{R}_{+},L^{\infty}(M))}$,
$||\lambda_u(t)||_{L^{\infty}(\mathbb{R}_{+})}$ and
$||\lambda_v(t)||_{L^{\infty}(\mathbb{R}_{+})}$ are not less than
$C$, where $C$ is only depending on upper bound of
$||g_u||_{H^1(M)}, ||g_v||_{H^1(M)}$, $||g_u||_{L^{\infty}(M)},
||g_v||_{L^{\infty}(M)}$. We calculate
\begin{eqnarray*}
\frac{1}{2}\frac{d}{dt}\int_{M}(u-v)^2dx
&=&\int_{M}(u-v)(u_t-v_t)dx\\
&=&\int_{M}(u-v)(\Delta(u-v)+\lambda_u(t)u-\lambda_v(t)v-(u^p-v^p))dx\\
&\leq&-\int_{M}|\nabla(u-v)|^2dx+\int_{M}(u-v)(\lambda_u(t)u-\lambda_v(t)v)dx
\end{eqnarray*}
Note that
\begin{eqnarray*}
&&\int_{M}(u-v)(\lambda_u(t)u-\lambda_v(t)v)dx\\
&=&(\lambda_u(t)-\lambda_v(t))\int_{M}(u-v)udx+\lambda_v(t)\int_{M}(u-v)^2dx\\
&\leq&|\lambda_u(t)-\lambda_v(t)|(\int_{M}(u-v)^2dx)^{\frac{1}{2}}(\int_{M}u^2dx)^{\frac{1}{2}}
+\lambda_v(t)\int_{M}(u-v)^2dx\\
&\leq&C|\lambda_u(t)-\lambda_v(t)|(\int_{M}(u-v)^2dx)^{\frac{1}{2}}
+C\int_{M}(u-v)^2dx,
\end{eqnarray*}
and
\begin{eqnarray}\label{lambda1}
|\lambda_u(t)-\lambda_v(t)|
&=&|\int_{M}((|\nabla u|^2-|\nabla v|^2)+(u^{p+1}-v^{p+1}))dx|\nonumber\\
&\leq&\int_{M}|\nabla(u-v)|(|\nabla u|+|\nabla
v|)dx+|\int_{M}(u-v)(\frac{u^{p+1}-v^{p+1}}{u-v})dx|\nonumber\\
&\leq&(\int_{M}|\nabla(u-v)|^2dx)^{\frac{1}{2}}(\int_{M}(|\nabla
u|+|\nabla
v|)^2dx)^{\frac{1}{2}}\nonumber\\
&&+(\int_{M}(u-v)^2dx)^{\frac{1}{2}}(\int_{M}(\frac{u^{p+1}-v^{p+1}}{u-v})^2dx)^{\frac{1}{2}}\nonumber\\
&\leq&C(\int_{M}|\nabla(u-v)|^2dx)^{\frac{1}{2}}+C(\int_{M}(u-v)^2dx)^{\frac{1}{2}}.
\end{eqnarray}
We have
\begin{eqnarray*}
&&\frac{1}{2}\frac{d}{dt}\int_{M}(u-v)^2dx\\
&\leq&-\int_{M}|\nabla(u-v)|^2dx+C^2(\int_{M}|\nabla(u-v)|^2dx)^{\frac{1}{2}}(\int_{M}(u-v)^2dx)^{\frac{1}{2}}\\
&&+C^2\int_{M}(u-v)^2dx+C\int_{M}(u-v)^2dx\\
&\leq&-\frac{1}{2}\int_{M}|\nabla(u-v)|^2dx+(\frac{C^4}{2}+C^2+C)\int_{M}(u-v)^2dx.
\end{eqnarray*}
By Gronwall inequality, we have
$$
||u-v||_{L^2}^2\leq ||g_u-g_v||_{L^2}^2
\exp((\frac{C^4}{2}+C^2+C)t).
$$
Further more,
\begin{eqnarray*}
\frac{1}{2}\frac{d}{dt}\int_{M}|\nabla(u-v)|^2dx
&=&\int_{M}\nabla(u-v) \cdot \nabla(u-v)_tdx\\
&=&-\int_{M}\Delta(u-v) \cdot (u-v)_tdx\\
&=&-\int_{M}\Delta(u-v)\cdot(\Delta(u-v)+\lambda_u(t)u-\lambda_v(t)v-(u^p-v^p))dx\\
&=&-\int_{M}(\Delta(u-v))^2dx+\int_{M}\nabla(u-v)\cdot\nabla(\lambda_u(t)u-\lambda_v(t)v)dx\\
&&+\int_{M}\Delta(u-v)\cdot (u^p-v^p)dx.
\end{eqnarray*}
Note that
\begin{eqnarray*}
&&\int_{M}\Delta(u-v)\cdot (u^p-v^p)dx\\
&\leq&\frac{1}{2}\int_{M}(\Delta(u-v))^2dx+\frac{1}{2}\int_{M}(u^p-v^p)^2dx\\
&=&\frac{1}{2}\int_{M}(\Delta(u-v))^2dx+\frac{1}{2}\int_{M}(u-v)^2(\frac{u^p-v^p}{u-v})^2dx\\
&\leq&\frac{1}{2}\int_{M}(\Delta(u-v))^2dx+\frac{C}{2}\int_{M}(u-v)^2dx
\end{eqnarray*}

Likewise,
\begin{eqnarray*}
&&\int_{M}\nabla(u-v)\cdot\nabla(\lambda_u(t)u-\lambda_v(t)v)dx\\
&=&(\lambda_u(t)-\lambda_v(t))\int_{M}\nabla(u-v)\cdot \nabla udx+\lambda_v(t)\int_{M}|\nabla(u-v)|^2dx\\
&\leq&|\lambda_u(t)-\lambda_v(t)|(\int_{M}|\nabla(u-v)|^2dx)^{\frac{1}{2}}(\int_{M}|\nabla
u|^2dx)^{\frac{1}{2}}
+\lambda_v(t)\int_{M}|\nabla(u-v)|^2dx\\
&\leq&C|\lambda_u(t)-\lambda_v(t)|(\int_{M}|\nabla(u-v)|^2dx)^{\frac{1}{2}}
+C\int_{M}|\nabla(u-v)|^2dx\\
&\leq&C^2(\int_{M}|\nabla(u-v)|^2dx+(\int_{M}|\nabla(u-v)|^2dx)^{\frac{1}{2}}(\int_{M}(u-v)^2dx)^{\frac{1}{2}})\\
& &+C\int_{M}|\nabla(u-v)|^2dx\\
&\leq&(\frac{3}{2}C^2+C)\int_{M}|\nabla(u-v)|^2dx+\frac{C^2}{2}\int_{M}(u-v)^2dx,
\end{eqnarray*}
where the second inequality follows by (\ref{lambda1}). Then we have
\begin{eqnarray*}
\frac{1}{2}\frac{d}{dt}\int_{M}|\nabla(u-v)|^2dx &\leq&
(\frac{3}{2}C^2+C)\int_{M}|\nabla(u-v)|^2dx+\frac{C^2+C}{2}\int_{M}(u-v)^2dx\\
&\leq&C'\int_{M}|\nabla(u-v)|^2dx,
\end{eqnarray*}
where the second inequality follows by Poincare inequality and $C'$
is a constant depending on the constant of Poincare inequality and
$C$. By Gronwall inequality, we have
$$
||\nabla(u-v)||_{L^2}^2\leq ||\nabla(g_u-g_v)||_{L^2}^2 \exp(C't).
$$
$\Box$

Finally, we give the proofs of Corollary \ref{cor1} and Corollary
\ref{cor2} below.

\textbf{Proof of Corollary \ref{cor1}.}

\noindent Since
$$
\frac{1}{2}\frac{d}{dt}\int_{M}|\nabla
u|^2dx+\int_{M}(u_t)^2dx=\int_{M}u_tAdx,
$$
we have
$$
\int_{M}|\nabla u|^2dx-\int_{M}|\nabla
g|^2dx+2\int^t_0\int_{M}(u_t)^2dxdt=2\int^t_0\int_{M}u_tAdxdt,
$$
Since $A\in L^2(\mathbb{R}_+, L^2(M))$, we get
\begin{eqnarray}\label{energe1}
\int^t_0\int_{M}(u_t)^2dxdt \leq \int_{M}|\nabla
g|^2dx+\int^t_0\int_{M}A^2dxdt\leq C.
\end{eqnarray}
So $\int^{\infty}_s\int_{M}(u_t)^2dxdt\to 0$ as $s\to \infty$.

By the arguments in Theorem \ref{thm1}, we have
$\lambda(t)=\int_{M}(|\nabla u|^2-uA) dx$ is continuous, uniformly
bounded in $t\in [0,\infty)$. Moreover, $u\in
L^{\infty}(\mathbb{R}_{+},H^1(M))$. Then we can take a subsequence
$\{t_i\}$ with $t_i\to \infty$ such that $u_i(x)=u(x,t_i)$,
$\lambda(t_i)\to\lambda_{\infty}$. By (\ref{energe1}) and Theorem
\ref{thm1}, we have
\begin{equation*}
\left\{
\begin{array}{l}
         u_i\to u_{\infty}\ \ \ \ \text{in}\ L^2(M), \\
         u_i\rightharpoonup u_{\infty} \ \ \ \text{in}\ H^1(M),\\
         \partial_t u_i-(\lambda(t_i)-\lambda_{\infty})u_i\to 0 \ \ \ \text{in}\ L^2(M).\\
\end{array}
\right.
\end{equation*}
Since $\partial_t u_i-(\lambda(t_i)-\lambda_{\infty})u_i=\Delta
u_i+\lambda_{\infty}u_i+A$, $u_{\infty}\in H^1$ solves the equation
$\Delta u_{\infty}+\lambda_{\infty}u_{\infty}+A=0$ in $M$ and
$\int_{M}|u_{\infty}|^2dx=1$.
 $\Box$

\textbf{Proof of Corollary \ref{cor2}.}

\noindent Since
$$
\frac{1}{2}\frac{d}{dt}\int_{M}|\nabla
u|^2dx=-\int_{M}(u_t)^2dx-\frac{1}{p+1}\frac{d}{dt}\int_{M}u^{p+1}dx,
$$
we have
\begin{equation}\label{energe2}
\lambda(t)+2\int^t_0\int_{M}|u_t|^2dxdt=\int_{M}(|\nabla
g|^2+\frac{2}{p+1}g^{p+1})dx +\frac{p-1}{p+1}\int_{M}u^{p+1}dx.
\end{equation}
By the arguments in Theorem \ref{thm2}, we have
$\lambda(t)=\int_{M}(|\nabla u|^2+u^{p+1}) dx$ is continuous,
uniformly bounded in $t\in [0,\infty)$. Moreover, $u\in
L^{\infty}(\mathbb{R}_{+},H^1(M))$ and $u\in
L^{\infty}(\mathbb{R}_{+},L^{p+1}(M))$. Then we can take a
subsequence $\{t_i\}$ with $t_i\to \infty$ such that
$u_i(x)=u(x,t_i)$, $\lambda(t_i)\to\lambda_{\infty}$. By
(\ref{energe2}) and Theorem \ref{thm2}, we have
\begin{equation*}
\left\{
\begin{array}{l}
         u_i\to u_{\infty}\ \ \ \ \text{in}\ L^2(M), \\
         u_i\rightharpoonup u_{\infty} \ \ \ \text{in}\ H^1(M)\ \text{and}\ L^p(M),\\
         \partial_t u_i-(\lambda(t_i)-\lambda_{\infty})u_i\to 0 \ \ \ \text{in}\ L^2(M).\\
\end{array}
\right.
\end{equation*}
Since $\partial_t u_i-(\lambda(t_i)-\lambda_{\infty})u_i=\Delta
u_i+\lambda_{\infty}u_i-u_i^p$, $u_i\in H^1$ solves the equation
$\Delta u_{\infty}+\lambda_{\infty}u_{\infty}-u_{\infty}^p=0$ in $M$
and $\int_{M}|u_{\infty}|^2dx=1$.
 $\Box$

\section{Gradient estimates}\label{sect3}
This section is devoted to the proofs of Theorem \ref{thm3} and
Theorem \ref{thm4}. We show in this section that the Harnack
quantity trick introduced in the fundamental work of P.Li and
S.T.Yau. The trick is to find a suitable Harnack quantity and apply
the maximum principle in a nice way.

\textbf{Proof of Theorem \ref{thm3}}. Let $u>0$ be a smooth solution
to the heat equation on $M\times [0,T)$
\begin{equation}\label{heat}
(\partial_t-\Delta) u=\lambda(t)u+A(x,t).
\end{equation}
Set
$$
w=log u.
$$
Then we have
\begin{equation}\label{ediff}
(\partial_t-\Delta) w=|\nabla w|^2+(\lambda+u^{-1}A).
\end{equation}

Following Li-Yau \cite{SY} we let $F=t(|\nabla
w|^2+a(\lambda+u^{-1}A)-aw_t)$ (where $a>1$) be the Harnack quantity
for (\ref{heat}). Then we have $$ |\nabla
w|^2=\frac{F}{t}-a(\lambda+u^{-1}A)+aw_t,
$$
$$
\Delta w=w_t-|\nabla
w|^2-(\lambda+u^{-1}A)=-\frac{F}{at}-(1-\frac{1}{a})|\nabla w|^2.
$$
and
$$
w_t-\Delta w=|\nabla
w|^2+(\lambda+u^{-1}A)=\frac{F}{t}+(1-a)(\lambda+u^{-1}A)+aw_t.
$$
 Note that
$$
(\partial_t-\Delta) w_t=2\nabla w\nabla
w_t+\frac{d}{dt}(\lambda+u^{-1}A).
$$
Using the Bochner formula, we have
 $$ (\partial_t-\Delta) |\nabla
w|^2=2\nabla w\nabla w_t-[2|D^2w|^2+2(\nabla w,\nabla\Delta
w)+2Ric(\nabla w,\nabla w)],
$$
and using (\ref{ediff}) we get
$$
(\partial_t-\Delta) |\nabla w|^2=2\nabla w\nabla (w_t-\Delta
w)-[2|D^2w|^2+2Ric(\nabla w,\nabla w)],
$$
which can be rewritten as
$$
(\partial_t-\Delta) |\nabla w|^2=2\nabla w\nabla
[\frac{F}{t}+(1-a)(\lambda+u^{-1}A)+aw_t]-[2|D^2w|^2+2Ric(\nabla
w,\nabla w)].
$$
Then we have $$ (\partial_t-\Delta) (|\nabla w|^2-aw_t) =2\nabla
w\nabla [\frac{F}{t}+(1-a)(\lambda+u^{-1}A)]
$$
$$-[2|D^2w|^2+2Ric(\nabla w,\nabla
w)]-a\frac{d}{dt}(\lambda+u^{-1}A).
$$
Hence
\begin{eqnarray*}
&&(\partial_t-\Delta) (|\nabla w|^2-aw_t+a(\lambda+u^{-1}A))\\
&=&(\partial_t-\Delta) (|\nabla w|^2-aw_t)+a(\partial_t-\Delta)
(\lambda+u^{-1}A))\\
&=& 2\nabla w\nabla [\frac{F}{t}+(1-a)(\lambda+u^{-1}A)]
-[2|D^2w|^2+2Ric(\nabla w,\nabla
w)]-a\frac{d}{dt}(\lambda+u^{-1}A)\\
& &+a(\partial_t-\Delta) (\lambda+u^{-1}A))\\
&=&2\nabla w\nabla [\frac{F}{t}+(1-a)(u^{-1}A)]
-[2|D^2w|^2+2Ric(\nabla w,\nabla w)]-a\Delta (u^{-1}A).
\end{eqnarray*}
Then we have
$$
(\partial_t-\Delta)F =\frac{F}{t}+2t\nabla w\nabla
[\frac{F}{t}+(1-a)(u^{-1}A)]
$$
$$-t[2|D^2w|^2+2Ric(\nabla w,\nabla
w)]-at\Delta (u^{-1}A).
$$
Assume that $$ \sup_{M\times[0,T]} F>0.
$$
Applying the maximum principle at the maximum point $(z,s)$, we then
have
$$ (\partial_t-\Delta)F\geq 0, \; \nabla F=0.
$$
In the following our computation is always at the point $(z,s)$. So
we get
\begin{equation}\label{key1}
\frac{F}{s}+2(1-a)s\nabla w\nabla (u^{-1}A)-s[2|D^2w|^2+2Ric(\nabla
w,\nabla w)]-as\Delta (u^{-1}A)\geq 0.
\end{equation}
That is
\begin{equation}\label{key2}
F-as^2\Delta (u^{-1}A)\geq 2(a-1)s^2\nabla w\nabla (u^{-1}A)+
s^2[2|D^2w|^2+2Ric(\nabla w,\nabla w)].
\end{equation}
Set
$$
\mu=\frac{|\nabla w|^2}{F}{(z,s)}.
$$
Then at $(z,s)$,$$ |\nabla w|^2=\mu F.
$$
Hence
$$
\nabla \frac{A}{u}=\frac{\nabla A}{u}-\frac{A\nabla
u}{u^2}=\frac{\nabla A}{u}-\frac{A}{u}\nabla w.
$$
So
$$
\nabla w \cdot \nabla \frac{A}{u}=\frac{\nabla w \cdot\nabla
A}{u}-\frac{A}{u}|\nabla w|^2\geq -\frac{|\nabla w||\nabla
A|}{u}-\frac{A}{u}|\nabla w|^2
$$
\begin{equation}\label{simplify1}
=-\frac{|\nabla A|}{u}\sqrt{\mu F}-\frac{A}{u}\mu F \geq
-\frac{1}{2}\frac{|\nabla A|^2}{u}-(\frac{1}{2}+A)\frac{\mu F}{u}.
\end{equation}
Further more, we have
\begin{eqnarray*}
(\partial_t-\Delta)(u^{-1}A)&=&\frac{1}{u}(\partial_t-\Delta)A-\frac{A}{u^2}(\partial_t-\Delta)u+
\frac{2}{u^2}\nabla u \cdot \nabla A-2\frac{A}{u^3}|\nabla u|^2\\
&=&\frac{1}{u}(\partial_t-\Delta)A-\frac{A}{u^2}(\lambda u+A)+
\frac{2}{u}\nabla w \cdot \nabla A-2\frac{A}{u}|\nabla w|^2\\
&\leq&\frac{1}{u}(\partial_t-\Delta)A-\frac{A}{u^2}(\lambda u+A)+
\frac{2}{u}\sqrt{\mu F}|\nabla A|-2\frac{A}{u}\mu F\\
&\leq&\frac{1}{u}(\partial_t-\Delta)A-\frac{A}{u}(\lambda+u^{-1}A)+
\frac{\mu F}{u}+\frac{|\nabla A|^2}{u}-2\frac{A}{u}\mu F,
\end{eqnarray*}
and
\begin{eqnarray*}
\partial_t(u^{-1}A)&=&\frac{A_t}{u}-\frac{A}{u^2}u_t\\
&=&\frac{A_t}{u}-\frac{A}{u}w_t\\
&=&\frac{A_t}{u}-\frac{A}{u}(\frac{1}{a}(|\nabla
w|^2-\frac{F}{s})+(\lambda+u^{-1}A))\\
&=&\frac{A_t}{u}-\frac{A}{u}\cdot
\frac{F}{a}(\mu-\frac{1}{s})-\frac{A}{u}(\lambda+u^{-1}A).
\end{eqnarray*}
Hence
\begin{equation}\label{simplify2}
-\Delta(u^{-1}A)=(\partial_t-\Delta)(u^{-1}A)-\partial_t(u^{-1}A)
\end{equation}
$$
\leq-\frac{\Delta A}{u}+\frac{|\nabla A|^2}{u}+\frac{1-2A}{u}\mu
F+\frac{A}{u}\cdot \frac{F}{a}(\mu-\frac{1}{s})
$$
$$
<-\frac{\Delta A}{u}+\frac{|\nabla A|^2}{u}+\frac{1-2A}{u}\mu
F+\frac{A}{u}\cdot \frac{F}{a}\mu. \ \ \ \ \ \ \ \ \
$$
Note that
$$
|D^2w|^2+Ric(\nabla w,\nabla w)\geq \frac{1}{n}|\Delta w|^2-K|\nabla
w|^2.
$$
So
$$
|D^2w|^2+Ric(\nabla w,\nabla w)\geq
\frac{1}{n}(\frac{F}{as}+(1-\frac{1}{a})|\nabla w|^2)^2-K|\nabla
w|^2
$$
\begin{equation}\label{simplify3}
=\frac{F^2}{n}(\frac{1}{as}+(1-\frac{1}{a})\mu)^2-K\mu F.
\end{equation}
Substituting (\ref{simplify1}) (\ref{simplify2}) and
(\ref{simplify3}) into (\ref{key2}), we get
$$
F+\frac{as^2}{u}(-\Delta A+|\nabla A|^2)+\mu F
\frac{s^2}{u}(a+(1-2a)A)
$$
$$
\geq -s^2(a-1)\frac{|\nabla A|^2}{u}-(a-1)s^2\frac{1+2A}{u}\mu F
+\frac{2F^2}{n}(\frac{1}{a}+(1-\frac{1}{a})\mu s)^2-2s^2K\mu F.
$$
Assume that
$$
F\geq \frac{as^2}{u}(-\Delta A+|\nabla A|^2)+ s^2(a-1)\frac{|\nabla
A|^2}{u},
$$
for otherwise we are done. Then we have
$$
2F+\mu F \frac{s^2}{u}(a+(1-2a)A)+(a-1)s^2\frac{1+2A}{u}\mu
F+2s^2K\mu F
$$
$$
\geq\frac{2F^2}{n}(\frac{1}{a}+(1-\frac{1}{a})\mu s)^2.
$$

Simplifying this inequality, we get
$$
\frac{2F}{n}\frac{1}{a^2}\leq \frac{2}{(1+(a-1)\mu s)^2}+\frac{\mu
s}{(1+(a-1)\mu s)^2}
$$
$$
\cdot s(u^{-1}(a+(1-2a)A)+u^{-1}(a-1)(1+2A)+2K).
$$
Hence we have the estimate for $F$ at $(z,s)$ such that
$$
F(z,s)\leq C(u^{-1},|A|,|\nabla A|,|\Delta A|,K,a,T),
$$
which is the desired gradient estimate. $\Box$

The idea proof of Theorem \ref{thm4} is similar to Theorem
\ref{thm3}.

\textbf{Proof of Theorem \ref{thm4}}. Let $u>0$ be a smooth solution
to the heat equation on $M\times [0,T)$
\begin{equation}\label{heat2}
(\partial_t-\Delta) u=\lambda(t)u-u^p.
\end{equation}
Set
$$
w=log u.
$$
Then we have
\begin{equation}\label{ediff2}
(\partial_t-\Delta) w=|\nabla w|^2+(\lambda-u^{p-1}).
\end{equation}

Following Li-Yau \cite{SY} we let $F=t(|\nabla
w|^2+a(\lambda-u^{p-1})-aw_t)$ (where $a>1$) be the Harnack quantity
for (\ref{heat2}). Then we have $$ |\nabla
w|^2=\frac{F}{t}-a(\lambda-u^{p-1})+aw_t,
$$
$$
\Delta w=w_t-|\nabla
w|^2-(\lambda-u^{p-1})=-\frac{F}{at}-(1-\frac{1}{a})|\nabla w|^2.
$$
and
$$
w_t-\Delta w=|\nabla
w|^2+(\lambda-u^{p-1})=\frac{F}{t}+(1-a)(\lambda-u^{p-1})+aw_t.
$$
 Note that
$$
(\partial_t-\Delta) w_t=2\nabla w\nabla
w_t+\frac{d}{dt}(\lambda-u^{p-1}).
$$
Using the Bochner formula, we have
 $$ (\partial_t-\Delta) |\nabla
w|^2=2\nabla w\nabla w_t-[2|D^2w|^2+2(\nabla w,\nabla\Delta
w)+2Ric(\nabla w,\nabla w)],
$$
and using (\ref{ediff2}) we get
$$
(\partial_t-\Delta) |\nabla w|^2=2\nabla w\nabla (w_t-\Delta
w)-[2|D^2w|^2+2Ric(\nabla w,\nabla w)],
$$
which can be rewritten as
$$
(\partial_t-\Delta) |\nabla w|^2=2\nabla w\nabla
[\frac{F}{t}+(1-a)(\lambda-u^{p-1})+aw_t]-[2|D^2w|^2+2Ric(\nabla
w,\nabla w)].
$$
Then we have $$ (\partial_t-\Delta) (|\nabla w|^2-aw_t) =2\nabla
w\nabla [\frac{F}{t}+(1-a)(\lambda-u^{p-1})]
$$
$$-[2|D^2w|^2+2Ric(\nabla w,\nabla
w)]-a\frac{d}{dt}(\lambda-u^{p-1}).
$$
Hence
\begin{eqnarray*}
&&(\partial_t-\Delta) (|\nabla w|^2-aw_t+a(\lambda-u^{p-1}))\\
&=&(\partial_t-\Delta) (|\nabla w|^2-aw_t)+a(\partial_t-\Delta)
(\lambda-u^{p-1}))\\
&=& 2\nabla w\nabla [\frac{F}{t}+(1-a)(\lambda-u^{p-1})]
-[2|D^2w|^2+2Ric(\nabla w,\nabla
w)]\\
& &-a\frac{d}{dt}(\lambda-u^{p-1})+a(\partial_t-\Delta) (\lambda-u^{p-1}))\\
&=&2\nabla w\nabla [\frac{F}{t}+(a-1)u^{p-1}]
-[2|D^2w|^2+2Ric(\nabla w,\nabla w)]+a\Delta u^{p-1}.
\end{eqnarray*}
Then we have
$$
(\partial_t-\Delta)F =\frac{F}{t}+2t\nabla w\nabla
[\frac{F}{t}+(a-1)u^{p-1}]
$$
$$-t[2|D^2w|^2+2Ric(\nabla w,\nabla
w)]+at\Delta u^{p-1}.
$$
Assume that $$ \sup_{M\times[0,T]} F>0.
$$
Applying the maximum principle at the maximum point $(z,s)$, we then
have
$$ (\partial_t-\Delta)F\geq 0, \; \nabla F=0.
$$
In the following our computation is always at the point $(z,s)$. So
we get
\begin{equation}\label{key3}
\frac{F}{s}+2(a-1)s\nabla w\nabla u^{p-1}-s[2|D^2w|^2+2Ric(\nabla
w,\nabla w)]+as\Delta u^{p-1}\geq 0.
\end{equation}
That is
\begin{equation}\label{key4}
F+2(a-1)s^2\nabla w\nabla u^{p-1}+as^2\Delta u^{p-1}\geq
s^2[2|D^2w|^2+2Ric(\nabla w,\nabla w)].
\end{equation}
Set
$$
\mu=\frac{|\nabla w|^2}{F}{(z,s)}.
$$
Then at $(z,s)$,$$ |\nabla w|^2=\mu F.
$$
Hence
\begin{equation}\label{simplify4}
\nabla w\nabla u^{p-1}=(p-1)u^{p-1}|\nabla w|^2= (p-1)u^{p-1}\mu F.
\end{equation}
Since
$$
\frac{\Delta
u}{u}=\frac{u_t}{u}-(\lambda-u^{p-1})=w_t-(\lambda-u^{p-1})=\frac{1}{a}(-\frac{F}{s}+|\nabla
w|^2).
$$
We have
\begin{eqnarray}\label{simplify5}
\Delta u^{p-1} &=& \nonumber (p-1)u^{p-1}\frac{\Delta
u}{u}+(p-1)(p-2)u^{p-3}|\nabla u|^2\\\nonumber
&=&(p-1)u^{p-1}\frac{1}{a}(-\frac{F}{s}+|\nabla
w|^2)+(p-1)(p-2)u^{p-1}|\nabla w|^2\noindent\\
&=&(p-1)u^{p-1}\frac{1}{a}(-\frac{F}{s}+\mu F)+(p-1)(p-2)u^{p-1}\mu
F.
\end{eqnarray}
Note that
$$
|D^2w|^2+Ric(\nabla w,\nabla w)\geq \frac{1}{n}|\Delta w|^2-K|\nabla
w|^2.
$$
So
$$
|D^2w|^2+Ric(\nabla w,\nabla w)\geq
\frac{1}{n}(\frac{F}{as}+(1-\frac{1}{a})|\nabla w|^2)^2-K|\nabla
w|^2
$$
\begin{equation}\label{simplify6}
=\frac{F^2}{n}(\frac{1}{as}+(1-\frac{1}{a})\mu)^2-K\mu F.
\end{equation}
Substituting (\ref{simplify4}) (\ref{simplify5}) and
(\ref{simplify6}) into (\ref{key4}), we get
$$
F+2s^2(a-1)(p-1)u^{p-1}\mu F+s(p-1)u^{p-1}F(\mu s -1)
$$
$$
+as^2(p-1)(p-2)u^{p-1}\mu F \geq
\frac{2F^2}{n}(\frac{1}{a}+(1-\frac{1}{a})\mu s)^2-2s^2K\mu F.
$$
Simplifying this inequality, we get
$$
\frac{2F}{n}\frac{1}{a^2}\leq \frac{1}{(1+(a-1)\mu s)^2}+\frac{\mu
s}{(1+(a-1)\mu s)^2}
$$
$$
\cdot s(2(a-1)(p-1)u^{p-1}+(p-1)u^{p-1}+a(p-1)|p-2|u^{p-1}+2K)
$$
Hence we have the estimate for $F$ at $(z,s)$ such that
$$
F(z,s)\leq C(u^{p-1},K,a,p,T),
$$
which is the desired gradient estimate. $\Box$

\end{document}